\documentclass[12pt]{article}
\usepackage{amsmath,amsfonts,latexsym,amsthm,amssymb}
\usepackage{bbm}
\newcommand{\D}{D^\alpha}

\numberwithin{equation}{section}

\DeclareMathOperator{\const}{const}
\DeclareMathOperator{\dist}{dist}

\DeclareMathOperator{\supp}{supp}
\begin{document}

\newtheorem{lem}{Lemma}
\newtheorem{teo}{Theorem}
\newtheorem{cor}{Corollary}
\newtheorem{prop}{Proposition}

\numberwithin{teo}{section}
\numberwithin{lem}{section}
\numberwithin{cor}{section}
\numberwithin{prop}{section}

\pagestyle{plain}
\title{The  Vladimirov-Taibleson Operator: Inequalities, Dirichlet Problem, Boundary H\"older Regularity}
\author{Anatoly N. Kochubei
\\ \footnotesize Institute of Mathematics,\\
\footnotesize National Academy of Sciences of Ukraine,\\
\footnotesize Tereshchenkivska 3, Kyiv, 01024 Ukraine\\
\footnotesize E-mail: \ kochubei@imath.kiev.ua}
\date{}
\maketitle

\vspace*{3cm}
\begin{abstract}
We study the Vladimirov-Taibleson operator, a model example of a pseudo-differential operator acting on real- or complex-valued functions defined on a non-Archimedean local field. We prove analogs of classical inequalities for fractional Laplacian, study the counterpart of the Dirichlet problem including the property of boundary H\"older regularity of solutions.
\end{abstract}

\vspace{2cm}


{\bf Acknowledgement:} This work was funded in part under the State Budget of Ukraine, Program 6541230 "Support to Priority Trends of Scientific Research", and also under a fellowship of the Universities for Ukraine (U4U) program.

\newpage
\section{Introduction}

The  Vladimirov-Taibleson operator $\D$ of $p$-adic fractional differentiation is a model example of a pseudo-differential operator acting on real- or complex-valued functions defined on a non-Archimedean local field $K$, its subsets or, for the multi-dimensional case, on the space $K^n$.

A number of important results on spectral properties of $\D$, its perturbations and generalizations, as well as the theory of related partial (pseudo-) differential equations are covered by the monographs \cite{AKS,K2001,VVZ,T,Z,KKZ} and many recent papers, such as \cite{BGPW,BCW,BGHH,K2014,K2020,K2021,Z2022,TZ} and others.

To some extent, the theory of this operator is parallel to that of the fractional Laplacian of real analysis. While the topological and geometric properties of a non-Archimedean local field are quite different from those of $\mathbb R$, this parallelism is a reflection of deep unity of mathematics where similar properties of quite different objects happen to be found as consequences of parallel algebraic structures.

The simplest and most important example of a non-Archimedean local field is the field of $p$-adic numbers, appearing, by Ostrowski's theorem, as the only possible alternative to $\mathbb R$ as a completion of the field of rational numbers. Therefore the non-Archimedean mathematics is one of the main branches of mathematics as a whole, now appearing in a variety of applications.

In this paper, we prove non-Archimedean analogs of some classical inequalities including the weighted positivity property \cite{Ma} and various kinds of fractional Sobolev and Poincar\'e inequalities (see, for example, \cite{BV,Ma,RO}). As in the classical situation \cite{RO}, the Poincar\'e inequality is a basic tool for studying the Dirichlet problem, formulated appropriately for nonlocal operators. On the other hand, non-Archimedean analogs of classical inequalities are interesting in themselves. Note that the non-Archimedean potential theory was initiated by Haran \cite{H}.

The investigation of boundary H\"older regularity of solution is a natural task prompted by comparisons with the classical theory of elliptic equations \cite{GT,LZLH}. For the non-Archimedean case, there is also a different motivation: the evident examples of open sets, like balls and spheres, are simultaneously closed, thus not possessing a boundary. The nontrivial examples are either punctured disks or infinite disjoint unions of clopen (= closed and open) sets \cite{K2001}. While it is obvious what domans in $\mathbb R^n$ are good (this is defined in terms of smoothness of the boundary), properties of open sets in the non-Archimedean case are formulated in different geometric terms. Here we follow, with necessary modifications, the geometric method suggested by Lian et al \cite{LZLH}.

The structure of this paper is as follows. In Section 2, we collect necessary preliminaries about local fields $K$, the structure of open sets in $K^n$, the representation of $K^n$ in terms of the unramified extension of $K$ \cite{Ta,K2001,K2021}, Sobolev spaces of complex-valued functions on local fields \cite{GKR,GK1,GK2}. In Section 3, we prove analogs of classical inequalities. In Section 4, we study the Dirichlet problem for the operator $\D$, prove the existence of its weak solutions and the comparison theorem, a substitute of the maximum principle. Section 5 is devoted to the boundary H\"older regularity of solutions.

\section{Preliminaries}

{\bf 2.1. Local fields.} A non-Archimedean local field is a non-discrete totally disconnected locally compact topological field. Such a field $K$ is isomorphic either to a finite extension of the field $\mathbb Q_p$ of $p$-adic numbers (here $p$ is a prime number), if $K$ has characteristic zero, or to the field of formal Laurent series with coefficients from a finite field, if $\operatorname{char}K>0$. For basic notions and results regarding local fields see, for example, \cite{Se1,We,K2001}. We consider only non-Archimedean local fields.

Any local field is endowed with an absolute value $|\cdot |_K$, such that:
1) $|x|_K=0$ if and only if $x=0$, 2) $|xy|_K=|x|_K\cdot |y|_K$,
3) $|x+y|_K\le \max (|x|_K,|y|_K)$. The last property called the ultrametric one implies that $|x+y|_K=|x|_K$, if $|y|_K<|x|_K$.

The ring $O=\{ x\in K:\ |x|_K\le 1\}$ is called the ring of integers of $K$. For $K=\mathbb Q_p$, we write $\mathbb Z_p$ instead of $O$. The ideal $P=\{ x\in K:\ |x|_K<1\}$ contains such an element $\beta$ that $P=\beta O$. The quotient ring $\bar{K}=O/P$ is
a finite field called the residue field. The absolute value is called normalized, if $|\beta |_K=q^{-1}$ where $q$ is the cardinality of $O/P$. Unless stated otherwise, the absolute values used below are assumed normalized. Such absolute values take the values $q^N$, $N\in \mathbb Z$. In the case $K=\mathbb Q_p$, the field of $p$-adic numbers, $\beta =p$ (where $p$ is seen as an element) and $q=p$ (as a natural number).

The additive group of a local field $K$ is self-dual, so that the Fourier analysis on $K$ is similar to the classical one. Let $\chi$ be a fixed non-constant additive character on
$K$, which is assumed having rank zero, so that $\chi (x)\equiv 1$ for $x\in O$, while $\chi
(x_0)\ne 1$ for some $x_0\in K$ with $|x_0|_K=q$.

The Fourier transform of a complex-valued function $f\in L^1(K)$ is defined as
$$
(\mathcal F)(\xi )=\widehat{f}(\xi )=\int\limits_K\chi (x\xi )f(x)\,dx,\quad \xi
\in K,
$$
where $dx$ is the Haar measure on the
additive group of $K$ normalized in such a way that the measure
of $O$ equals 1.
If $\mathcal F f=\widehat{f}\in L^1(K)$, then the inversion rule
$$
f(x)=\int\limits_K\chi (-x\xi )\widehat{f}(\xi )\,d\xi ,
$$
is valid. We will denote $\widetilde{f}=\mathcal F^{-1}f$.

The Fourier transform preserves the Bruhat-Schwartz space $\mathcal D(K)$ of test functions, consisting of locally constant functions with compact supports. The local constancy of a function $f: K\to \mathbb C$ means the existence of such an integer $k$ that for any $x\in K$
$$
f(x+x')=f(x), \quad \text{whenever $|x'|\le q^{-k}$}.
$$

The above Fourier analysis is extended easily to functions on $K^n$. The natural non-Archimedean norm on $K^n$ is
$$
|(x_1,\ldots ,x_n)|_{K^n}=\max\limits_{1\le j\le n}|x_j|_K.
$$
The Fourier transform extends to the dual space $\mathcal D'(K)$ (and to its multi-dimensional counterpart $\mathcal D'(K^n)$) called the space of Bruhat-Schwartz distributions.

\medskip
{\bf 2.2. Field extensions and their use in analysis.}  If a local field $K$ is a subfield of a local field $L$, then $L$ is called an extension of $K$ (which is denoted $L/K$). Consider $L$ as a vector space over $K$. An extension $L/K$ is called finite, if the space $L$ is finite-dimensional over $K$. Its dimension is called the degree of the extension.

An operator of multiplication in a finite extension $L$ by an element $\xi$ can be considered as a linear operator in the $K$-vector space, so that $\operatorname{Tr}(\xi)$ is defined. The extension is called separable, if the linear function $\xi \mapsto \operatorname{Tr}(\xi)$ does not vanish identically. All finite extensions of a field of characteristic zero are separable. The above notion of separability makes sense also for finite fields $\bar{K},\bar{L}$.

A finite extension $L/K$ is called unramified, if $\bar{L}/\bar{K}$ is a separable extension of the same degree as $L/K$. Any local field $K$ has a unique (up to isomorphism) unramified extension of any given degree $n\ge 1$. Any prime element $\beta$ of the field $K$ is also a prime element of any unramified extension. If $L$ is an unramified extension of $K$ of a degree $n$, then the cardinality of the residue field equals $q^n$ where $q$ is the cardinality of the residue field of $K$.

As a vector space over $K$, the unramified extension $L$ of degree $n$ has a canonical basis consisting of representatives of a basis in $\bar{L}$ over $\bar{K}$. If $x\in L$ has the coefficients $x_1,\ldots ,x_n\in K$ of the expansion with respect to the canonical basis, then the normalized absolute value $|x|_L$ has the representation \cite{Ta,K2021}
\begin{equation}
\label{2.1}
|x|_L=\left( \max\limits_{1\le j\le n}|x_j|_K\right)^n.
\end{equation}

An automorphism $\sigma$ of the field $L$ is called an automorphism of the extension $L/K$, if $\sigma (a)=a$ for all $a\in K$. A finite extension $L/K$ is called a Galois extension, if the order of its group of automorphisms $G$ coincides with the degree of the extension. In this case $G$ is called the Galois group of the extension.

In the important case of an unramified extension $L/K$, the group $G$ is cyclic. Its generator $F$ is called the Frobenius automorphism. By the construction of the absolute value on $L$ (\cite{Se1}, Chapter II, \S2) and the expression of the norm map in terms of the Galois group (\cite{Se1}, Chapter V, \S2), together with the fact that the prime element in $K$ remains prime in $L$, we see that $F$ preserves the absolute value on $L$.

Let $L$ be an unramified extension of degree $n$ of a local field $K$. Taking into account (\ref{2.1}), we see the expansion with respect to a canonical basis in $L$ defines an isometric linear isomorphism between $L$ and $K^n$. In various applications (see, for example, \cite{K2021}), it is convenient to reduce problems for multi-dimensional operators acting on functions $K^n\to \mathbb C$, to one-dimensional operators on functions $L\to \mathbb C$ where $L/K$ is an unramified extension of degree $n$.

\medskip
{\bf 2.3. Structure of open sets.} Let $\Omega\subset K$ be an open subset of a local field $K$. The set $\Omega$ can be represented as a union
\begin{equation}
\label{2.2}
\bigcup_{k=1}^N V_k,\quad N\le \infty,
\end{equation}
of non-intersecting balls $V_k=B(x_k,r_k)=\{ x\in K: |x-x_k|_K\le q^{r_k}\},r_k\in \mathbb Z$. The set (\ref{2.2}) is compact, if and only if $N<\infty$. See \cite{Se} for further investigation of this case.

Turning to the non-compact case, we note that any local field is a separable metric space (\cite{Sch}, Exercise 19.B (iii)), and any separable metric space has the Lindel\"of property -- its every covering has a countable subcovering. Therefore we may assume that the disjoint union in (\ref{2.2}) is countable. It is known (Proposition 3.1 in \cite{K2001}) that an open set $\Omega$ is closed, if and only if the sequence $\{ x_k\}$ has no finite limit points. Now we prove a more precise result.

\medskip
\begin{prop}
The boundary $\partial \Omega=\operatorname{closure}(\Omega)\setminus \Omega$ coincides with the set of all finite limit points of the sequence of centers $\{x_k\}_1^\infty$.
\end{prop}

\medskip
{\it Proof}. Note first that a limit point of the sequence of centers never belongs to $\Omega$. Indeed, otherwise it would belong to one of the balls, and then an infinite subsequence of the centers would belong to that ball. However the latter property contradicts the disjointness of the balls. Thus, the above set of limit points belongs to $\partial \Omega$.

Conversely, let $x_0\in \partial \Omega$. There exists a sequence $\{ y_j\}\subset \Omega$, $y_j\to x_0$. Taking a subsequence if necessary, we may assume that $|y_j-x_{k_j}|_K\le q^{r_{k_j}}$ where all the numbers $k_j$ are different.

In addition, we have $r_{k_n}\to -\infty$, as $n\to \infty$. Indeed, otherwise there exists a subsequence $\{r_j'\} \subset \{ r_{k_n}\}$, such that $r_j'\ge R>-\infty$. For the corresponding subsequences $\{y_j'\}$ and $\{x_j'\}$, we get
$$
\left| y_{j_1}'-y_{j_2}'\right|_K=\left| (y_{j_1}'-x_{j_1}')+(x_{j_1}'-x_{j_2}')+(x_{j_2}'-y_{j_2}')\right|_K=\left| x_{j_1}'-x_{j_2}'\right|_K \ge \max\{ q^{r_{j_1}'},q^{r_{j_2}'}\}\ge q^R,
$$
which contradicts the convergence of the sequence $\{ y_j\}$.

Since $r_{k_n}\to -\infty$, for any $\varepsilon >0$, there exists such a number $j_1$ that $\left| x_{k_j}-y_j\right|_K<\varepsilon$ for $j\ge j_1$. On the other hand, there exists such a number $j_2$ that $\left| y_j-x_0\right|_K<\varepsilon$ for $j\ge j_2$. Now
$$
\left| x_{k_j}-x_0\right|_K\le \max\{ \left| x_{k_j}-y_j\right|_K,\left| y_j-x_0\right|_K\}<\varepsilon,
$$
as $j\ge \max\{j_1,j_2\}$. This means that $x_0$ is a limit point of a sequence $\{ x_n\}$.$\quad \blacksquare$

\medskip
The ``textbook examples'' of non-Archimedean open sets, like balls and spheres, are clopen and have no boundaries. The simplest nontrivial example is a punctured unit ball $\mathbb Z_p\setminus \{ 0\}$, for which the decomposition (\ref{2.2}) has the explicit form
$$
\mathbb Z_p\setminus \{ 0\}=\bigcup_{n=0}^\infty \left( \bigcup_{k=1}^{p-1}\left(kp^n+p^{n+1}\mathbb Z_p\right) \right)
$$
(see Remark 4 in \cite{LLF}). Here the centers of the balls are the points $kp^n$, and the sequence of centers tends to 0.

\medskip
{\bf 2.4. The Vladimirov operator.} On a test function $\varphi \in \mathcal D(K)$, the fractional differentiation operator $\D$, $\alpha >0$ is defined as follows:
\begin{equation}
\label{2.3}
(\D \varphi )(x)=\mathcal F^{-1}\left[ |\xi |_K^\alpha (\mathcal F(\varphi))(\xi )\right] (x), \quad x\in K.
\end{equation}
The operator $\D$ admits a hypersingular integral representation
\begin{equation}
\label{2.4}
\left( D^\alpha \varphi \right) (x)=\int\limits_K \mathcal K(y)[\varphi (x)-\varphi (x+y)]\,dy,
\end{equation}
where
$$
\mathcal K(y)=\frac{q^\alpha -1}{1-q^{-\alpha -1}}|y|_K^{-\alpha-1}.
$$
The expression (\ref{2.4}) makes sense for wider classes of functions.

For the multi-dimensional case, the natural generalization (often called the Taibleson operator) is defined as a pseudo-differential operator $D^\alpha_{K^n}$ with the symbol
$$
|(\xi_1, \ldots ,\xi_n)|_{K^n}^\alpha,\quad |(\xi_1, \ldots ,\xi_n)|_{K^n}=\max\limits_{1\le j\le n}|\xi_j|_K.
$$
In this case, we have a hypersingular integral representation similar to (\ref{2.4}), with the integration over $K^n$ and
$$
\mathcal K(y)=\frac{q^\alpha -1}{1-q^{-\alpha -n}}|y|_K^{-\alpha-n}.
$$
This operator can be interpreted as the one-dimensional Vladimirov operator $D^{\alpha/n}$ over the unramified extension of degree $n$ of $K$; see \cite{K2021}.

\medskip
{\bf 2.5. Sobolev spaces.} The Sobolev type spaces related to the operator $\D$ were introduced by Taibleson \cite{T}; see also \cite{Kim}. A more general case of locally compact abelian groups was studied by G\'orka et al \cite{GKR,GK1,GK2}. Here we will not use other Sobolev-type spaces introduced by Z\'u\~niga-Galindo; see \cite{KKZ}.

The Sobolev space $H^\alpha (K^n)$ where $K$ is a local field, $\alpha >0$, consists of such $f\in L^2 (K^n)$ that
$$
\|f\|_{H^\alpha (K^n)}=\left\{ \int\limits_{\xi \in K^n}|\widehat{f}(\xi )|^2(1+|\xi|^2_{K^N})^\alpha d\xi \right\}^{1/2}<\infty
$$
where $\widehat{f}=\mathcal Ff$.

The imbedding
$$
H^\alpha (K^n)\hookrightarrow L^{\alpha^*}(K^n),\quad \alpha^*=\frac{2\gamma}{\gamma -\alpha}\  \left( \gamma >\max \left( \frac{n}2,\alpha \right) \right),
$$
holds in the local field situation (\cite{GKR}, Theorem 5). This means that
$$
\|u\|_{L^{\alpha^*}(K^n)}\le C\|u\|_{H^\alpha (K^n)} \text{ for all $u\in H^\alpha (K^n)$}
$$
(here and below $C$ denotes various positive constants).

In the ``fractional'' case, that is $0<\alpha <1$, there is an equivalent Aronszain-Gagliardo-Slobodecki norm
$$
\|u\|^2_{AGS}=\|u\|^2_{L^2(K^n)}+[u]^2_\alpha
$$
where
$$
[u]_\alpha^2=\int\limits_{K^n}\int\limits_{K^n}\frac{|u(x)-u(y)|^2}{|x-y|_{K^n}^{2\alpha +n}}dx\,dy,
$$
and we have a more refined estimate called the fractional Sobolev inequality.

\medskip
\begin{prop}
For every $u\in H^\alpha (K^n)$, $0<\alpha <\min(1,\frac{n}2)$,
\begin{equation}
\label{2.5}
\|u\|_{L^{\alpha^*}(K^n)}\le C[u]_\alpha,\quad \alpha^*=\frac{2n}{n-2\alpha},
\end{equation}
where $C$ does not depend on $u$.
\end{prop}

\medskip
The {\it proof} for the non-Archimedean case is a slight modification of Brezis' proof for $\mathbb R^n$ (see \cite{BV,Po}). Let $u\in H^\alpha (K^n)$. For each $x,y\in K^n$,
$$
|u(x)|\le |u(x)-u(y)|+|u(y)|,
$$
so that we get integrating in $y\in B(x,q^l)$ that
$$
q^{ln}|u(x)|\le \int\limits_{|y-x|_{K^n}\le q^l} |u(x)-u(y)|\,dy+\int\limits_{|y-x|_{K^n}\le q^l} |u(y)|\,dy,
$$
so that
$$
|u(x)|\le q^{-ln}\int\limits_{|y-x|_{K^n}\le q^l} |u(x)-u(y)|\,dy+q^{-ln}\int\limits_{|y-x|_{K^n}\le q^l} |u(y)|\,dy,
$$

By the H\"older inequality,
\begin{multline*}
q^{-ln}\int\limits_{|y-x|_{K^n}\le q^l} |u(x)-u(y)|\,dy\le q^{-ln/2}\left\{ \int\limits_{|y-x|_{K^n}\le q^l} |u(x)-u(y)|^2\,dy\right\}^{1/2} \\
\le q^{\alpha l}\left\{ \int\limits_{|y-x|_{K^n}\le q^l} \frac{|u(x)-u(y)|^2}{|x-y|_{K^n}^{2\alpha +n}}\,dy\right\}^{1/2}.
\end{multline*}
The H\"older inequality yields also the estimate
$$
q^{-ln}\int\limits_{|y-x|_{K^n}\le q^l} |u(y)|\,dy \le \left\{ q^{-ln}\int\limits_{|y-x|_{K^n}\le q^l} |u(y)|^r\,dy\right\} ^{1/r},
$$
for any $r\in [1,\infty )$. Therefore
\begin{equation}
\label{2.6}
|u(x)|\le q^{\alpha l}\left\{ \int\limits_{K^n} \frac{|u(x)-u(y)|^2}{|x-y|_{K^n}^{2\alpha +n}}\,dy\right\}^{1/2} +q^{-ln/r}\left\{\int\limits_{K^n}|u(y)|^r\,dy\right\}^{1/r},
\end{equation}
for any $l\in \mathbb Z, r\in [1,\infty )$.

We may assume that $u$ is bounded. The general case is then considered \cite{Po} using the standard truncation argument. For a bounded $u$, the right-hand side of (\ref{2.6}) is finite for $r>2$, for almost all $x$.

For an arbitrary $z>0$, there exists such $l\in \mathbb Z$ that $q^l\le z\le q^{l+1}$. Then $q^{l\alpha}\le z^\alpha$,
$$
q^{-ln/r}=q^{-(l+1)n/r}\cdot q^{n/r}\le q^{n/r}z ^{-n/r},
$$
and it follows from (\ref{2.6}) that
\begin{equation}
\label{2.7}
|u(x)|\le C\left( az^\alpha +bz^{-n/r}\right)
\end{equation}
for any $z\ge 0$. Here
$$
a=\left\{ \int\limits_{K^n} \frac{|u(x)-u(y)|^2}{|x-y|_{K^n}^{2\alpha +n}}\,dy\right\}^{1/2},\quad b=\left\{\int\limits_{K^n}|u(y)|^r\,dy\right\}^{1/r},
$$

Minimizing the right-hand side of (\ref{2.7}) as it was done in \cite{BV,Po} and integrating we come to the inequality (\ref{2.5}). $\qquad \blacksquare$

\medskip
The above fractional Sobolev space is a special case of the Besov spaces on local fields studied in \cite{Ka}.

\medskip
\section{Inequalities}

{\bf 3.1. The Poincar\'e inequality.} Let $\Omega \subset K^n$ be a bounded open set. Denote
 \begin{equation}
\label{3.1}
X=\left\{ H^\alpha (K^n):\ u\equiv 0 \text{ in $K^n\setminus \Omega$} \right\} ,\quad 0<\alpha <1.
\end{equation}

The next result is similar to the Poincar\'e type inequality for nonlocal operators of real analysis; see the inequality (3.4) in \cite{RO}.

\medskip
\begin{teo}
There exists such a positive constant $C$ that for any $u\in X$,
\begin{equation}
\label{3.2}
\int\limits_\Omega |u(x)|^2\,dx \le C\int\limits_{K^n}\int\limits_{K^n}\frac{|u(x)-u(x+y)|^2}{|y|_{K^n}^{n+2\alpha}}\,dx\,dy.
\end{equation}
\end{teo}

\medskip
{\it Proof.} Let us apply the H\"older inequality to the integral on the left in (\ref{3.2}). Since the set $\Omega$ is bounded, we have
$$
\int\limits_\Omega |u(x)|^2\,dx \le C\left[ \int\limits_{K^n}(|u(x)|^2)^{\frac{n}{n-2\alpha}}dx\right]^{\frac{n-2\alpha}n}.
$$
By the inequality (\ref{2.5}),
$$
\left\{ \int\limits_{K^n}(|u(x)|^2)^{\frac{n}{n-2\alpha}}dx\right\}^{\frac{n-2\alpha}n}\le C\int\limits_{K^n}\int\limits_{K^n}\frac{|u(x)-u(x+y)|^2}{|y|_{K^n}^{n+2\alpha}}\,dx\,dy,
$$
and we obtain (\ref{3.2}). $\qquad \blacksquare$

\medskip
{\bf 3.2. The fractional Poincar\'e-Wirtinger inequality} (compare with \cite{TJZ}). Let $B_N^{(n)}=\left\{ x\in K^n:\ |x|_{K^n}\le q^N\right\}$, $N\in \mathbb Z$. Consider the Sobolev space $H^\alpha (B_N^{(n)})$ with the norm \linebreak $\|u\|^2_{\alpha,N,n}=\|u\|^2_{L^2(B_N^{(n)}})+[u]^2_{\alpha,N,n}$ where $0<\alpha <1$,
$$
[u]_{\alpha,N,n}^2=\int\limits_{B_N^{(n)}}\int\limits_{B_N^{(n)}}\frac{|u(x)-u(y)|^2}{|x-y|_{K^n}^{2\alpha +n}}\,dx\,dy.
$$

We begin with the case $n=1$, and in this case we drop $n=1$ from the notations.

\medskip
\begin{lem}
Let $V$ be a closed subspace in $H^\alpha (B_N)$, which does not contain nonzero constants. Then for any $u\in V$,
\begin{equation}
\label{3.3}
\|u\|_{L^2(B_N)}\le C(\alpha,N)[u]_{\alpha,N}.
\end{equation}
\end{lem}

\medskip
{\it Proof.} Suppose the opposite. Then there exists a sequence $\{ u_m\}\subset V$, such that \linebreak $\|u_m\|_{L^2(B_N)}> m[u_m]_{\alpha,N}$ for each $m\ge 1$. Let $v_m=\|u_m\|_{L^2(B_N)}^{-1}u_m$. Then $\|v_m\|_{L^2(B_N)}=1$,
$[v_m]_{\alpha,N}<\dfrac1m$. Taking a subsequence if necessary, we may assume that the sequence $\{ v_m\}$ is weakly convergent in $H^\alpha (B_N)$.

Below we will use the identity
\begin{equation}
\label{3.4}
|x|_K^\alpha =\frac{1-q^\alpha}{1-q^{-\alpha -1}}\int\limits_K|\xi |_K^{-\alpha -1}[\chi (x\xi)-1]\,d\xi.
\end{equation}

To prove (\ref{3.4}), we use the Riesz kernel \cite{K2001,VVZ}
$$
f_\gamma (x)=\frac{|x|_K^{\gamma -1}}{\Gamma_K(\gamma)},\quad x\in K,\gamma >0, \gamma \ne 1,
$$
where $\Gamma_K(\gamma )=\dfrac{1-q^{\gamma -1}}{1-q^{-\gamma}}$. Considering $f_\gamma$ as a distribution from $\mathcal D'(K)$ we have the identity for its Fourier
transform, $\widetilde{f_\gamma}(\xi )=|\xi|_K^{-\gamma}$. Now, for any $\varphi \in \mathcal D(K)$ we denote $\psi =\mathcal F^{-1}\varphi$ and find that
\begin{multline*}
\langle |x|_K^\alpha ,\varphi \rangle =\langle f_{\alpha +1},\varphi\rangle \Gamma_K (\alpha+1)=\langle \widetilde{f_{\alpha +1}},\psi \rangle \Gamma_K (\alpha+1)=\Gamma_K (\alpha+1)\langle |x|_K^{-\alpha -1},\psi (x)-\psi (0)\rangle \\
=\Gamma_K (\alpha+1)\int\limits_K  |x|_K^{-\alpha -1}\,dx\int\limits_K[\chi (-x\xi )-1]\varphi (\xi )\,d\xi.
\end{multline*}
Applying the Fubini theorem we come to (\ref{3.4}).

The ball $B_N$ is an additive locally compact Abelian group, and the space $H^\alpha (B_N)$ can be interpreted in terms of the Pontryagin duality. While in the proof of (\ref{3.4}) we used the harmonic analysis on $K$, now we switch to harmonic analysis on $B_N$ preserving the notations for additive characters and the Fourier transform.

The dual group $\widehat{B}_N$ to $B_N$ is isomorphic to the discrete group $K/B_{-N}$ consisting of the cosets
$$
\xi +B_{-N}=\beta^m\left( r_0+r_1\beta +\cdots +r_{N-m-1}\right) +B_{-N},\quad m\in \mathbb Z,m<N,
$$
where $r_j$ belongs to a complete set of representatives in $O$ of the elements of the residue field $O/P$. On $\widehat{B}_N$, there is a normalized discrete measure $d\xi$ satisfying the Plancherel identity (see, for example, \cite{K2018}).

It is shown in \cite{GK2} (the formula (12)) that
\begin{equation}
\label{3.5}
[u]^2_{\alpha,N}=\int\limits_{\widehat{B}_N}|\widehat{u}(\xi +B_{-N})|^2A(\xi )\,d\xi
\end{equation}
where
\begin{equation}
\label{3.6}
A(\xi )=\int\limits_{B_N}\frac{|\chi (z\xi)-1|^2}{|z|_K^{1+2\alpha}}\,dz.
\end{equation}

By (\ref{3.4}), $A(\xi )=\const \cdot |\xi |_K^\alpha$, and it follows from (\ref{3.5}) that the above definition of the norm in $H^\alpha (B_N)$ is equivalent to the definition of the Sobolev norm in terms of the Fourier transform on $B_N$. On the basis of the latter definition, it is proved in \cite{GKR} (Theorem 11) that the imbedding $H^\alpha (B_N)\hookrightarrow L^2(B_N)$ is compact.

Therefore the sequence $\{v_m\}$, weakly convergent in $H^\alpha (B_N)$, converges strongly in $L^2(B_N)$ to a certain function $v$, and $\|v\|_{L^2(B_N)}=1$. On the other hand, $[\cdot ]_{\alpha, N}$ is a norm on $V$. Since a norm is lower semicontinuous (see \cite{Br}, page 61), $[v]_{\alpha,N}\le \liminf [v_m]_{\alpha,N}=0$, that is $[v]_{\alpha,N}=0$, so that $v=\const$, and by our assumption, $v\equiv 0$, which contradicts the equality $\|v\|_{L^2(B_N)}=1$. $\qquad \blacksquare$

The analogue of the Poincar\'e-Wirtinger inequality is as follows.

\begin{teo}
For any $u\in H^\alpha (B_N)$,
\begin{equation}
\label{3.7}
\|u-\bar{u}\|_{L^2(B_N)}\le C[u]_{\alpha,N},\quad\bar{u}=q^{-N}\int\limits_{B_N}u(x)\,dx.
\end{equation}
\end{teo}

\medskip
{\it Proof}. Let us consider the subspace
$$
V=\left\{ u\in H^\alpha (B_N):\ \int\limits_{B_N}u(x)\,dx=0\right\}.
$$
Then $u-\bar{u}\in V$, $[u-\bar{u}]_{\alpha,N}=[u]_{\alpha,N}$. Substituting into (\ref{3.3}) we obtain (\ref{3.7}). $\qquad \blacksquare$

\medskip
Turning to the multi-dimensional case where we consider functions on $K^n$, we identify $K^n$ with the unramified extension $L$ of the field $K$ of degree $n$. Then by (\ref{2.1}),
$$
B_N^{(n)}=\{ x\in K^n:\ |x|_{K^n}\le q^N\}=\{ x\in L:\ |x|_L^{1/n}\le q^N\} =\{ x\in L:\ |x|_L\le r^N\}
$$
where $r$ is the cardinality of the residue field of $L$.

Next, consider the seminorm $[u]_{L,\frac{\alpha}n,N}$ of the form
$$
[u]_{L,\frac{\alpha}n,N}^2=\int\limits_{|x|_L\le r^N}\int\limits_{|y|_L\le r^N} \frac{|u(x)-u(y)|^2}{|x-y|^{1+2\alpha/n}}\,dx\,dy.
$$
Using (\ref{2.1}) once more we see that
$$
[u]_{L,\frac{\alpha}n,N}=[u]_{\alpha,n,N}.
$$
Applying Theorem 3.2, we obtain the following result.

\medskip
\begin{cor}
For any $u\in H^\alpha (B_N^{(n)})$,
$$
\|u-\bar{u}\|_{L^2(B_N^{(n)})}\le C[u]_{\alpha,N,n},\quad\bar{u}=q^{-Nn}\int\limits_{B_N^{(n)}}u(x)\,dx.
$$
\end{cor}

\medskip

Note that a little weaker inequality was proved by a different method in \cite{BGHH} for a more general framework of ultrametric spaces.

\medskip
{\bf 3.3. The Poincar\'e inequality with $\alpha$-harmonic capacity} (compare with Section 6.5 of \cite{Ma}).

Let $K$ be a local field. As before, we identify $K^n$ with the unramified extension $L$ of the field $K$ of degree $n$. Let $r=q^n$ is the cardinality of the residue field of $L$. A prime element $\beta$ of $K$ is prime also for $L$. For any element $x\in L$, $|x|_L=r^N$, $N\in \mathbb Z$, we can write the canonical representation
$$
x=\beta^{-N}(\xi_1+\xi_2\beta +\cdots )
$$
where $\xi_j$ belong to a complete set of representatives of residue classes from $O_L/P_L$.

Using the same notations as in the previous section, we consider the ball $B_N\subset L$, $B_N=\{ x\in L: |x|_L\le r^N\}$. Then the inclusion $x\in B_N$ means that
\begin{equation}
\label{3.8}
x-\beta^{-N}\xi_1=\beta^{-N+1}(\xi_2+\xi_3\beta+\cdots ).
\end{equation}
The ball $B_N$ is represented as a disjoint union of $r$ balls $B_{N-1}(\xi_1)$ of the radius $r^{N-1}$ described by (\ref{3.8}) with fixed $\xi_1$. The Frobenius automorphism $F$ transposes these balls.

Let $f$ be a Lipschitz function on the ball $B_{N-1}(\xi_1^0)$ where $\xi_1^0$ is an arbitrary fixed element. For any $\xi_1\ne \xi_1^0$,, there exists $\nu \in \{ 1,\ldots ,r-1\}$, such that $F^{-\nu} (\xi_1^0)=\xi_1$, so that
$F^{-\nu}:\  B_{N-1}(\xi_1^0)\to B_{N-1}(\xi_1)$, and the function
$$
f_\nu (x)=f(F^\nu (x )),\quad x\in B_{N-1}(\xi_1),
$$
is defined.

Since the Frobenius automorphism and its powers preserve the absolute value, and the distance between points of different balls equals $r^N$, the functions $f_\nu$ define a Lipschitz function $f^*$ on $B_N$. The mapping $f\mapsto f^*$ is a continuous mapping of the spaces of Lipschitz functions.

If $e$ is a compact subset of $B_{N-1}(\xi_1^0)$, such that $\dist (\supp f,e)>0$, then $\dist (\supp f^*,e)>0$.

The $\alpha$-{\it capacity} $\operatorname{cap}_\alpha (e,B_N)$ is defined as
$$
\operatorname{cap}_\alpha (e,B_N)=\inf \left\{ [u]_{\alpha,N}^2:\ u\in \mathcal D(B_N),u=1 \text{ in a neighborhood of $e$}\right\}.
$$

\medskip
\begin{teo}
Let $e$ be a compact subset of $B_N$. For any real-valued function $u\in \mathcal D(B_{N-1}(\xi_1^0))$, such that $\dist (\supp u,e)>0$, we have the inequality
\begin{equation}
\label{3.9}
\operatorname{cap}_\alpha (e,B_N)\|u\|^2_{L^2(B_{N-1}(\xi_1^0))}\le C[u]_{\alpha,B_{N-1}(\xi_1^0)}.
\end{equation}
\end{teo}

\medskip
{\it Proof}. Let $f=1-u$, $f^*$ be the above extension onto $B_N$. Suppose that $\eta \in \mathcal D(B_N)$, $\eta=1$ on a neighborhood of $B_{N-1}(\xi_1^0)$. Then
\begin{equation}
\label{3.10}
\operatorname{cap}_\alpha (e,B_N)\le \| \eta f^*\|^2_{\alpha,B_{N-1}(\xi_1^0)}\le C \| 1-u\|^2_{\alpha,B_{N-1}(\xi_1^0)}
\end{equation}
where $C$ does not depend on $u$. Here we used the invariance of the absolute value and the Haar measure with respect to the Frobenius automorphism. The double integral over
$$
F^{\nu_1}(B_{N-1}(\xi_1^0))\times F^{\nu_2}(B_{N-1}(\xi_1^0)),\quad \nu_1\ne \nu_2,
$$
is estimated via the $L^2$-norm.

It follows from (\ref{3.10}) that
\begin{equation}
\label{3.11}
\operatorname{cap}_\alpha (e,B_N)\le C\inf\left\{  \| 1-u\|^2_{\alpha,B_{N-1}(\xi_1^0)}:\ u\in \mathcal D(B_{N-1}(\xi_1^0)),\ \dist (\supp u,e)>0\right\}.
\end{equation}

Let
$$
M=\left\{ \frac1{\operatorname{mes}B_{N-1}(\xi_1^0)}\int\limits_{B_{N-1}(\xi_1^0)}u^2(x)\,dx\right\}^{1/2}
$$
where $\operatorname{mes}B_{N-1}(\xi_1^0)=r^{n-1}$. By (\ref{3.11}),
$$
\operatorname{cap}_\alpha (e,B_N)\le C\|1-M^{-1}u\|_{\alpha,B_{N-1}(\xi_1^0)}=CM^{-2}[u]_{\alpha,B_{N-1}(\xi_1^0)}
+C\|1-M^{-1}u\|_{L^2(B_{N-1}(\xi_1^0))},
$$
that is
\begin{equation}
\label{3.12}
M^2\operatorname{cap}_\alpha (e,B_N)\le C[u]_{\alpha,B_{N-1}(\xi_1^0)}
+C\|M-u\|_{L^2(B_{N-1}(\xi_1^0))}.
\end{equation}

Denote
$$
\bar{u}_1=\left( \operatorname{mes}B_{N-1}(\xi_1^0)\right)^{-1}\int\limits_{B_{N-1}(\xi_1^0)}u(x)\,dx.
$$
We may assume that $\bar{u}_1\ge 0$. Then by the Cauchy inequality, $\bar{u}_1\le M$,
$$
\|\bar{u}_1\|_{L^2(B_{N-1}(\xi_1^0))}=\bar{u}_1( \operatorname{mes}B_{N-1}(\xi_1^0))^{1/2},
$$
$$
\|u\|_{L^2(B_{N-1}(\xi_1^0))}-\|\bar{u}_1\|_{L^2(B_{N-1}(\xi_1^0))}=( \operatorname{mes}B_{N-1}(\xi_1^0))^{1/2}(M-\bar{u}_1),
$$
so that
$$
( \operatorname{mes}B_{N-1}(\xi_1^0))^{1/2}(M-\bar{u}_1)\le \|u-\bar{u}_1\|_{L^2(B_{N-1}(\xi_1^0))},
$$
and we obtain from Theorem 3.2 that
\begin{multline*}
\|M-u\|_{L^2(B_{N-1}(\xi_1^0))}\le \|M-\bar{u}_1\|_{L^2(B_{N-1}(\xi_1^0))}+\|u-\bar{u}_1\|_{L^2(B_{N-1}(\xi_1^0))}\le 2\|u-\bar{u}_1\|_{L^2(B_{N-1}(\xi_1^0))}\\
\le C[u]_{\alpha,B_{N-1}(\xi_1^0)}.
\end{multline*}
Now the inequality (\ref{3.12}) implies (\ref{3.9}). $\qquad \blacksquare$

\medskip
{\bf 3.4. Weighted positivity} -- compare with Section 8.3 in \cite{Ma}.

In this section we prove an integral identity for the operator $\D$ implying, in particular, its positivity in a weighted Hilbert space whose weight is the fundamental solution for $\D$.

\medskip
\begin{lem}
If $u,v$ are real-valued functions belonging to $\mathcal D(K^n)$, $0<\alpha <n$, then
\begin{multline}
\label{3.13}
u(x)(\D v)(x)+v(x)(\D u)(x)-(\D (u\cdot v))(x)\\
=a_\alpha\int\limits_{K^n}\frac{[u(x)-u(x+y)][v(x)-v(x+y)]}{|y|_{K^n}^{n+\alpha}}\,dy,\quad x\in K^n,
\end{multline}
where $a_\alpha =\dfrac{q^{\alpha n}-1}{1-q^{-\alpha-n}}$.
\end{lem}

{\it Proof}. Let us multiply the expression in the right-hand side of (\ref{3.13}) by a function $\varphi \in \mathcal D(K^n)$ and integrate over $K^n$. We use the well-known properties of the Fourier transform (see \cite{VVZ}): if $\varphi,\psi\in  \mathcal D(K^n)$, then
$$
\int \varphi\overline{\psi}\,dx=\int \widehat{\varphi}\overline{\widehat{\psi}}\,dx;\quad \widehat{\varphi *\psi}=\widehat{\varphi}\widehat{\psi};\quad \widehat{\varphi \psi}=\widehat{\varphi}*\widehat{\psi}.
$$
We find that
\begin{multline}
\label{3.14}
\int [u(x)-u(x+y)][v(x)-v(x+y)]\varphi (x)\,dx\\
=\int \overline{\mathcal F_x[u(x)-u(x+y)](\xi )}\mathcal F_x\{ [v(x)-v(x+y)]\varphi (x)\}(\xi )\,d\xi .
\end{multline}

Next,
$$
\mathcal F_x[u(x)-u(x+y)](\xi )=\int \chi (x\cdot \xi)[u(x)-u(x+y)]\,dx=[1-\chi (-y\cdot \xi)]\widehat{u}(\xi).
$$
The right-hand side of (\ref{3.14}) equals
\begin{multline}
\label{3.15}
\int [1-\chi (y\cdot \xi)]\overline{\widehat{u}(\xi)}\mathcal F_x\{ [v(x)-v(x+y)]*\varphi (x)\}(\xi )\,d\xi  \\
=\int\int [1-\chi (-y\cdot \xi)][1-\chi (-y\cdot \eta)]\overline{\widehat{u}(\xi)}\widehat{v}(\eta)\widehat{\varphi}(\xi -\eta)\,d\xi\,d\eta.
\end{multline}

The integral in $y$ is evaluated using the identity
\begin{equation}
\label{3.16}
|x|_{K^n}^\alpha =\frac{1-q^\alpha}{1-q^{-\alpha -n}}\int\limits_{K^n}|\xi|_{K^n}^{-\alpha -n}[\chi (x\cdot \xi)-1]\,d\xi,
\end{equation}
a multi-dimensional version of the identity (\ref{3.4}). In fact, (\ref{3.16}) can be obtained from (\ref{3.4}) using the approach based on the unramified extension $L$ of degree $n$ of the field $K$. As in Section 3.2, we consider the operator $D^{\alpha/n}$ on $L$ and use the relations $r=q^n$, $|x|_L=|x|^n_{K^n}$ (see (\ref{2.1})).

Using (\ref{3.16}), we get
\begin{multline*}
\int [1-\chi (-y\cdot \xi)][1-\chi (-y\cdot \eta)]|y|_{K^n}^{-\alpha -n}\,dy\\
= -\int {(\chi (y\cdot \xi)-1)+(\chi (-y\cdot \eta)-1)-(\chi (y\cdot (\xi -\eta))-1)}|y|_{K^n}^{-\alpha -n}\,dy \\
=\frac1{a_\alpha}\left( |\xi|^\alpha_{K^n}+|\eta|^\alpha_{K^n}-|\xi -\eta|^\alpha_{K^n}\right).
\end{multline*}

Taking into account (\ref{3.14}) and (\ref{3.15}) we obtain that the result of multiplying in (\ref{3.13}) by $\varphi$, with subsequent integration, is the expression
$$
\int\int \left( |\xi|^\alpha_{K^n}+|\eta|^\alpha_{K^n}-|\xi -\eta|^\alpha_{K^n}\right) \overline{\widehat{u}(\xi)}\widehat{v}(\eta)\widehat{\varphi}(\xi -\eta)\,d\xi\,d\eta.
$$
Here $|\eta|^\alpha_{K^n}\widehat{v}(\eta)=(\mathcal F(\D v))(\eta)$,
$$
\int (\mathcal F(\D v))(\eta)\widehat{\varphi}(\xi -\eta)\,d\eta = ((\mathcal F(\D v))*\widehat{\varphi}) (\xi)=(\mathcal F(\D v\cdot \varphi))(\xi),
$$
so that
\begin{multline*}
\int\int |\eta|^\alpha_{K^n} \overline{\widehat{u}(\xi)}\widehat{v}(\eta)\widehat{\varphi}(\xi -\eta)\,d\xi\,d\eta =\int \overline{\widehat{u}(\xi)}(\mathcal F(\D v\cdot \varphi))(\xi)
=\int u(x)(\D v)(x)\varphi (x)\,dx.
\end{multline*}
Similarly,
$$
\int\int |\xi|^\alpha_{K^n} \overline{\widehat{u}(\xi)}\widehat{v}(\eta)\widehat{\varphi}(\xi -\eta)\,d\xi\,d\eta=\int v(x)(\D u)(x)\varphi (x)\,dx.
$$
Finally,
\begin{multline*}
\int\int |\xi -\eta|^\alpha_{K^n} \overline{\widehat{u}(\xi)}\widehat{v}(\eta)\widehat{\varphi}(\xi -\eta)\,d\xi\,d\eta =\int\int  |\tau|^\alpha_{K^n} \overline{\widehat{u}(\xi)}\widehat{v}(\xi -\tau)\widehat{\varphi}(\tau)\,d\xi\,d\tau \\
=\int (\mathcal F(\D \varphi))(\tau)(\overline{\widehat{u}}*\overline{\widehat{v})(\tau)})\,d\tau=\int (\D \varphi)(x)u(x)v(x)\,dx=\int \varphi (x)(\D (uv))(x)(x)\,dx.
\end{multline*}
Since $\varphi$ is arbitrary, we come to (\ref{3.13}). $\qquad \blacksquare$

\medskip
Note that the Riesz kernel $f_{\alpha}$ (see the above Section 3.2 and Section VIII.4 in \cite{VVZ} where the multi-dimensional case is considered) is a fundamental solution for the operator $\D$, that is $\D f_\alpha=\delta$; this is
easily verified using the Fourier transform \cite{K2001,VVZ}. The explicit expression for $f_\alpha$ is
$$
f_\alpha (x)=\frac{|x|_{K^n}^{\alpha -1}}{\Gamma_{K^n}(\alpha)},\quad \Gamma_{K^n}(\alpha)=\frac{1-q^{\alpha -n}}{1-q^{-\alpha}}.
$$
It is important that $f_\alpha (x)\ge 0$, if $0<\alpha <n$.

Consider the equality (\ref{3.13}) for $u=v$, multiply both sides by $f_\alpha$ and integrate. We obtain the next result.

\medskip
\begin{teo}[weighted positivity]
Let $0<\alpha <n$. For every real-valued $u\in \mathcal D(K^n)$, the following equality is valid; both sides are nonnegative:
$$
2\int (\D u)(x)u(x)f_\alpha (x)\,dx=u(0)^2+a_\alpha \int\int \frac{|u(x)-u(y)|^2}{|x-y|_{K^n}^{n+\alpha}}f_\alpha (x)\,dx\,dy.
$$
\end{teo}

\medskip
\section{Dirichlet Problem}

{\bf 4.1. The case of a homogeneous boundary condition.} The Poincar\'e type inequality (\ref{3.2}) implies weak solvability of the Dirichlet problem

\begin{gather}
\D u=f \text{ in $\Omega \subset K^n$}, \\
u=0 \text{ in $K^n\setminus \Omega$},
\end{gather}
where $\Omega$ is a bounded open subset of $K^n$. As in analysis on $\mathbb R^n$ \cite{RO}, for a nonlocal operator, the boundary condition is set on the complement of $\Omega$.

Let $X_{K^n}$ be the space of functions $u(x),x\in K^n$, such that $u\equiv 0$ on $K^n\setminus \Omega$,
$$
\int\limits_{K^n}\int\limits_{K^n}\frac{|u(x)-u(x+y)|^2}{|y|^{n+2\alpha}_{K^n}}\,dx\,dy<\infty.
$$
$X_{K^n}$ is a Hilbert space with the inner product
$$
(v,w)=\frac12 \int\limits_{K^n}\int\limits_{K^n}(v(x)-v(x+y))(w(x)-w(x+y))\mathcal K(y)\,dx\,dy.
$$

The weak formulation of the problem (4.1)-(4.2) is as follows. Suppose $f\in X'_{K^n}$ (the dual space); we write our equation as
$$
(u,\varphi)=\int\limits_\Omega f(x)\varphi (x)\,dx \text{ for all $\varphi \in X_{K^n}$}.
$$

The existence of a unique weak solution is a consequence of the Lax-Milgram theorem (see Corollary 5.8 in \cite{Br}).

\medskip
{\bf 4.2. Inhomogeneous boundary condition.} The Dirichlet problem
\begin{gather}
\D u=0 \text{ in $\Omega \subset K^n$}, \\
u=g \text{ in $K^n\setminus \Omega$},
\end{gather}
has been studied by Haran \cite{H} as a part of his $p$-adic potential theory. Here we reproduce some results from \cite{H} we use in a sequel. As before, we extend easily some results from $K$ to $K^n$ using the unramified extension technique.

Let $\Omega=B_{-N-1}^{(n)}$, $K^n\setminus \Omega=\{ x\in K^N:\ |x|_{K^n}\ge q^{-N}\}$. Let $g$ be a continuous function on $K^n\setminus \Omega$, such that
$$
\int\limits_{K^n\setminus \Omega}|g(x)|\frac{dx}{|x|_{K^n}^{n+\alpha}}<\infty.
$$
Define the function $u$ as the extension of $g$ onto $K^n$ by the constant
$$
\frac{1-q^{-\alpha}}{1-q^{-n}}q^{-N\alpha}\int\limits_{K^n\setminus \Omega}g(x)\frac{dx}{|x|_{K^n}^{n+\alpha}}
$$
Then the function $u$ is a solution of (4.3)-(4.4) (\cite{H}, page 935).

In other words, the above expression is a counterpart of the classical Poisson kernel. It is used (\cite{H}, Section 7.7) to obtain an expression for the Green function $G_\Omega^\alpha (x,y)$ corresponding to the problem (4.1)-(4.2) with $\Omega=B_{-N-1}$.

For our purposes, it is sufficient to list some of its properties. Namely, for a ball $\Omega$,
\begin{equation}
\label{4.5}
0\le G_\Omega^\alpha (x,y)\le C|x-y|_{K^n}^{\alpha -n}
\end{equation}
(see \cite{H}, page 935), where $C$ does not depend on $\Omega$,
$G_\Omega^\alpha (x,y)=0$ for $x\in K^n\setminus \Omega$, $y\in \Omega$, and also for $y\in K^n\setminus \Omega$ and any $x$. We do not touch a more involved theory based on the notion of regular boundary points.

{\bf 4.3. Comparison theorem.} In the $p$-adic case, a comparison theorem is different from real counterparts where the technique is based typically on the non-existence of continuous step functions; see, for example, Lemma 9 in \cite{LL}. For the non-Archimedean local field $K$, the space $\mathcal D(K)$ of (continuous) step functions is dense in $L^2(K)$. Therefore non-Archimedean comparison theorems contain additional assumptions.

For linear equations, it suffices to prove positivity of solutions.

\medskip
\begin{teo}
Let $u$ be a continuous weak solution of the Dirichlet problem
\begin{gather}
\D u=f \text{ in $\Omega \subset K^n$}, \\
u=g \text{ in $K^n\setminus \Omega$},
\end{gather}
where $f$ and $g$ are continuous functions, $f\ge 0$ on $\Omega$, $g\ge 0$ on $\Omega^c=K^n\setminus \Omega$, $\Omega$ is a bounded open subset of $K^n$ satisfying at least one of the following conditions:

(i) $\Omega$ has a nonempty boundary $\partial\Omega$;

(ii) For each $x\in \Omega$,
$$
x+\Omega \subset \Omega;\quad x+\Omega^c\subset \Omega^c.
$$
\medskip
Then $u\ge 0$ on $\Omega$.
\end{teo}

\medskip
{\it Proof}. The function $u$ is a weak solution in the sense of  Section 4.1, satisfying the appropriate identity for an arbitrary test function $\varphi$.

In order to specify $\varphi$, write $u=u^+-u^-$ on $\Omega$, that is
$$
u^+=\max\{ u,0\}\mathbbm{1}_\Omega, u^-=\max\{ -u,0\}\mathbbm{1}_\Omega.
$$
where $\mathbbm{1}_\Omega$ is the indicator function of the set $\Omega$.
If $u$ is not nonnegative, then $u^-$ is not identically zero, and we may set $\varphi =u^-$.

We have
\begin{equation}
\label{4.8}
\iint\limits_{(K^n\times K^n)\setminus (\Omega^c\times \Omega^c)}(u(x)-u(z))(\varphi (x)-\varphi (z))\mathcal K(z-x)\,dx\,dz=\int\limits_\Omega f\varphi \,dx
\end{equation}
where the right-hand side is nonnegative. On the other hand, the left-hand side of (\ref{4.8}) equals
$$
\int\limits_\Omega \int\limits_\Omega (u(x)-u(z))(u^-(x)-u^-(z))\mathcal K(z-x)\,dx\,dz + 2\int\limits_\Omega dx\int\limits_{\Omega^c}(u(x)-g(z))u^-(x)\mathcal K(z-x)\,dz.
$$

Next, $u^+(x)u^-(x)=0$, so that
\begin{multline*}
(u(x)-u(z))(u^-(x)-u^-(z))=(u^+(x)-u^+(z))(u^-(x)-u^-(z))-(u^-(x)-u^-(z))^2\\
=-(u^+(x)u^-(z)+u^+(z)u^-(x))-(u^-(x)-u^-(z))^2 \le 0,
\end{multline*}
and returning to (\ref{4.8}) we see that the left-hand side is less than or equal to zero. To avoid the contradiction, we have to consider a possibility that $u^-(x)$ equals identically a nonnegative constant $C$, that is $u(x)\equiv -C$ on
$\Omega$

Under our assumption (i), there exists a point $x_0\in \partial \Omega$. It is known (\cite{Bo}, 1.1.6) that $x_0$ belongs also to the boundary of $\Omega^c$. Since $u$ is continuous on $K^n$, its value $u(x_0)$ must coincide with $\lim\limits_{\Omega \ni x\to x_0} u(x)$ and with $\lim\limits_{\Omega^c \ni x\to x_0} u(x)$. The first of these limits equals $-C\le 0$, the second one equals $g(x_0)\ge 0$. Therefore, $C=0$, $u^-=0$, and we have come to a contradiction.

Under the assumption (ii), for any $x\in \Omega$,
\begin{multline*}
(\D u)(x)=\int\limits_{\Omega}[u(x)-u(x+y)]\mathcal K(y)\,dy+\int\limits_{\Omega^c}[-C-u(x+y)]\mathcal K(y)\,dy\\
=-C\int\limits_{\Omega^c}\mathcal K(y)\,dy -\int\limits_{\Omega^c}g(x+y)\mathcal K(y)\,dy\le -C\int\limits_{\Omega^c}\mathcal K(y)\,dy <0.
\end{multline*}
This contradiction completes the proof. $\qquad\blacksquare$

\medskip
\section{Boundary Regularity}

In this section, we consider the Dirichlet problem (4.6)-(4.7) where $0<\alpha <1$, $\Omega$ is a bounded open set with nonempty boundary. To be definite, we assume that $0\in \partial \Omega$. To simplify notations, we consider the one-dimensional case, $\Omega \subset K$; the general situation can be studied, as before, using the unramified extension of degree $n$.

We also assume that $f\in L^s(\Omega)$, $s>\frac1\alpha$, $g$ is bounded and belongs to the H\"older space $C^\delta (0)$, $g(0)=0$.

Our task is to find conditions on $\Omega$, under which a continuous solution of the problem (4.6)-(4.7) satisfies the estimate
\begin{equation}
\label{5.1}
|u(x)|\le C|x|_{K}^\gamma,\quad \gamma >0,
\end{equation}
on a neighborhood of the origin $0\in \partial \Omega$.

Geometric properties of $\Omega$ are described by the following condition: there exists $\nu \in (0,1)$, such that for all $k>0$
\begin{equation}
\label{5.2}
\operatorname{mes}\left[(B(r_k)\setminus  B(r_{k+1}))\cap \Omega^c\right] \ge \nu r_k
\end{equation}
where $r_k=q^{-\lambda_k}$, $\{ \lambda_k\}$ is a sequence of natural numbers satisfying the quasi-geometric growth condition, that is
\begin{equation}
\label{5.3}
1<R_-\le \frac{\lambda_{k+1}}{\lambda_k}\le R_+<\infty,\quad k=1,2,\ldots;\ \lambda_1=1.
\end{equation}

\medskip
\begin{teo}
Let $u$ be a continuous solution of the Dirichlet problem (4.6)-(4.7). Denote
$$
M=\|u\|_{L^\infty (\Omega \cap B(1))}+ \|f\|_{L^s(\Omega\cap B(1))}+\|g\|_{C^\delta (0)}.
$$
Then the geometric condition (\ref{5.2}) implies the inequality
$$
|u(x)|\le \widehat{C}M|x|_K^\gamma,\text{ for some $\gamma >0$},
$$
where $\widehat{C}$ and $\gamma$ do not depend on $u$, valid in a neighborhood of the origin.
\end{teo}

\medskip{\it Proof}. Let us construct a sequence of nonnegative functions $\{ v_k\}_1^\infty $, such that $v_1\equiv M$, while for $k\ge 2$,
\begin{equation}
\label{5.4}
\left\{
\begin{array}{ccc}
\D v_k=|f| \text{ on $B(r_{(k-1)k_0})$};\\
-v_k\le u\le v_k\text{ on $K$};\\
v_k\le \widehat{C}Mq^{-\lambda_{kk_0}\rho} \text{ on $B(r_{kk_0})$}
\end{array}
\right.
\end{equation}
where $k_0\ge 1$ is a natural number, $0<\rho <\min (\alpha -\frac1s,\delta )$, $\widehat{C}$ does not depend on $u$ and $k$. Specific values of $k_0$ and $\rho$ (independent of $k$) will be chosen later.

For $k=2$, we set
\begin{equation}
\label{5.5}
v_2(x)=\int\limits_{B(r_{k_0})^c}P^\alpha_{r_{k_0}}(x,y)v_1(y)\,dy +\int\limits_{B(r_{k_0})}G^\alpha_{r_{k_0}}(x,y)|f|(y)\,dy
\end{equation}
where $P^\alpha_{r_{k_0}}(x,y)$ is the Poisson kernel for the ball $B(r_{k_0})$ understood in the sense of Section 4.2, $G^\alpha_{r_{k_0}}(x,y)$ is the Green function for the same ball.

Since the Poisson kernel, for any $x\in B(r_{k_0})$, defines a probability measure, the first summand in (\ref{5.5}) equals $M$ and can be extended continuously onto $K$ by this constant. By the ultrametric property and the Holder inequality, the second summand is less than or equal to $c_1Mq^{-\lambda_{k_0}(\alpha -\frac1s)}$ where $c_1>0$ does not depend on $k_0$. We have, for $0<\rho_1<\min (\alpha -\frac1s,\delta )$, that
$$
v_2(x)\le M+c_1Mq^{-\rho_1\lambda_{k_0}}=M+c_1Mr_{k_0}^{\rho_1}.
$$

By (\ref{5.3}),
\begin{gather*}
\frac{\lambda_{n+1}}{R_+}\le \lambda_n,\\
\frac{\lambda_{n+2}}{R_+}\le \lambda_{n+1},\\
\ldots \\
\frac{\lambda_{n+n}}{R_+}\le \lambda_{n+n-1}.
\end{gather*}
Multiplying these inequalities we find that
\begin{equation}
\label{5.6}
\frac{\lambda_{2n}}{R_+^n}\le \lambda_n,\quad n=1,2,\ldots .
\end{equation}

In particular, $\lambda_{k_0}\ge R_+^{-k_0}\lambda_{2k_0}$, so that
$$
q^{-\rho_1\lambda_{k_0}}\le q^{-\sigma \lambda_{2k_0}},\quad \sigma =\frac{\rho_1}{R_+^{k_0}}.
$$
Taking a larger constant $\widehat{C}$ we find that $v_2(x)$ satisfies the inequality from (\ref{5.4}) with $k=2$:
\begin{equation}
\label{5.7}
v_2(x)\le \widehat{C}Mq^{-\sigma \lambda_{2k_0}}, x\in B(r_{2k_0}).
\end{equation}

The inequality $-v_2\le u\le v_2$ is a consequence of Theorem 4.1.

Suppose that the estimate from (\ref{5.4}) is fulfilled for some $k\ge 2$. Let us prove it for $k+1$. Denote
$$
\mathfrak G(x)=\begin{cases}
Mr^\delta_{(k-1)k_0} & \text{ on $\left( B(r_{(k-1)k_0})\setminus B(r_{(kk_0})\right) \cap \Omega^c$};\\
 \widehat{C}Mr^\sigma_{(k-1)k_0} & \text{ on $\left( B(r_{(k-1)k_0})\setminus B(r_{kk_0})\right) \cap \Omega$};\\
 v_k & \text{ on $K\setminus \left( B(r_{(k-1)k_0})\setminus B(r_{kk_0})\right)$}.
 \end{cases}
 $$
 Then $-\mathfrak G\le u\le \mathfrak G$ on $B(r_{kk_0})^c$. Let
 \begin{equation}
\label{5.8}
v_{k+1}(x)=\int\limits_{B(r_{kk_0})^c}P^\alpha_{r_{kk_0}}(x,y)\mathfrak G(y)\,dy +\int\limits_{B(r_{kk_0})}G^\alpha_{r_{kk_0}}(x,y)|f|(y)\,dy,\quad x\in B(r_{kk_0}).
\end{equation}

Let us write the first summand in (\ref{5.8}) as $I_1+I_2$ where
$$
I_1=\int\limits_{B(r_{(k-1)k_0})\setminus B(r_{kk_0})}P^\alpha_{r_{kk_0}}(x,y)\mathfrak G(y)\,dy,\quad I_2=\int\limits_{B(r_{(k-1)k_0})^c}P^\alpha_{r_{kk_0}}(x,y)\mathfrak G(y)\,dy.
$$
Let
$$
A=(B(r_{(k-1)k_0})\setminus B(r_{kk_0}))\cap \Omega^c=\bigcup\limits_{j=(k-1)k_0}^{kk_0}(B(r_j)\setminus B(r_{j+1}))\cap \Omega^c.
$$
Using (\ref{5.2}) and the notation $\zeta (\alpha )=(1-q^{-\alpha})^{-1}$ we write
\begin{multline*}
\int\limits_A P^\alpha_{r_{kk_0}}(x,y)\,dy=\frac1{\zeta (\alpha)}q^{-\alpha \lambda_{kk_0}}\sum\limits_{j=(k-1)k_0}^{kk_0}\int\limits_{(B(r_j)\setminus B(r_{j+1}))\cap \Omega^c} \frac{dy}{|y|_K^{1+\alpha}}\\
=\frac1{\zeta (\alpha)}q^{-\alpha \lambda_{kk_0}}\sum\limits_{j=(k-1)k_0}^{kk_0}q^{(1+\alpha)\lambda_j}\int\limits_{(B(r_j)\setminus B(r_{j+1}))\cap \Omega^c} dy\\
\ge \frac1{\zeta (\alpha)}q^{-\alpha \lambda_{kk_0}}\nu \sum\limits_{j=(k-1)k_0}^{kk_0}q^{\alpha\lambda_j} \ge \frac{\nu}{\zeta (\alpha)}\overset{\text{def}}{=} \mu \in (0,1).
\end{multline*}

Let
 \begin{equation}
\label{5.9}
\mathfrak H=\frac{\mathfrak G-Mr_{(k-1)k_0}^\delta }{\widehat{C}Mr_{(k-1)k_0}^\sigma -Mr_{(k-1)k_0}^\delta}.
\end{equation}
Since $0<\sigma <\delta$, the denominator in (\ref{5.9}) is positive, provided $\widehat{C}\ge 1$. Therefore $\mathfrak H\le 0$ on $A$. Moreover, $\mathfrak H\le 1$ on $B(r_{(k-1)k_0})$. Therefore
\begin{multline*}
\int\limits_{B(r_{(k-1)k_0})\setminus B(r_{kk_0})} P^\alpha_{r_{kk_0}}(x,y)\mathfrak H(y)\,dy\\
=\int\limits_{B(r_{(k-1)k_0})\setminus B(r_{kk_0})} P^\alpha_{r_{kk_0}}(x,y)\,dy
+\int\limits_{B(r_{(k-1)k_0})\setminus B(r_{kk_0})} P^\alpha_{r_{kk_0}}(x,y)(\mathfrak H(y)-1)\,dy \\
\le \int\limits_{B(r_{kk_0})^c}P^\alpha_{r_{kk_0}}(x,y)\,dy - \int\limits_A P^\alpha_{r_{kk_0}}(x,y)\,dy\le 1-\mu.
\end{multline*}

Returning to $\mathfrak G$ and expressing $\mathfrak G$ via $\mathfrak H$, we see that
\begin{multline*}
I_1\le (1-\mu )\left[ \widehat{C}Mr_{(k-1)k_0}^\sigma -Mr_{(k-1)k_0}^\delta\right] +Mr_{(k-1)k_0}^\delta
\le [(1-\mu )\widehat{C}M+\mu M]r_{(k+1)k_0}^{\sigma/R_+^{k_0}}.
\end{multline*}

Let us estimate $I_2$. It follows by induction from (\ref{5.8}) and the properties of the Poisson kernel and Green function that $v_k(x)=M$ for $x\in B(r_{(k-1)k_0})^c$. Now
\begin{multline*}
I_2=M\int\limits_{B(r_{(k-1)k_0})^c}P^\alpha_{r_{kk_0}}(x,y)=\frac{\zeta (1)}{\zeta (\alpha)}Mq^{-\alpha \lambda_{kk_0}}\int\limits_{|y|_K\ge q^{-\lambda_{(k-1)k_0}+1}}\frac{dy}{|y|_K^{1+\alpha}}\\
=\frac{M}{\zeta (\alpha)}q^{-\alpha \lambda_{kk_0}}\sum\limits_{j=\lambda_{(k-1)k_0}+1}^\infty q^{-j\alpha} =Mq^{-\alpha \lambda_{kk_0}+\alpha \lambda_{(k-1)k_0}-1}.
\end{multline*}

Like in the previous estimates, we find from (\ref{5.3}) that
$$
\lambda_n\le \frac{\lambda_{n+l}}{R_-^l},\quad l\ge 0.
$$
In particular,
$$
\lambda_{(k-1)k_0}\le \frac{\lambda_{kk_0}}{R_-^{k_0}},
$$
so that
$$
I_2\le Mq^{-\alpha \lambda_{kk_0}(1-\frac1{R_-^{k_0}})R_-^{k_0}\lambda_{(k+1)k_0}}\le \widehat{C}Mq^{-\rho \lambda_{(k+1)k_0}},
$$
if $k_0$ is chosen in such a way that
$$
\rho=\alpha (1-R_-^{-k_0})R_-^{-k_0}<\min (\alpha -\frac1s,\delta).
$$

Next, using the H\"older inequality we get as before that
$$
\int\limits_{B(r_{kk_0})}G^\alpha_{r_{kk_0}}(x,y)|f|(y)\,dy\le c_1Mr_{kk_0}^{\alpha -\frac1s}\le c_1Mr_{kk_0}^\rho\le  c_1Mr_{(k+1)k_0}^{\rho (R_-^{k_0})^{-1}}.
$$
Taking together several latest estimates we can choose $k_0$ and $\rho$ (not depending on $k$) in such a way that the inequality in (\ref{5.4}) holds for the step $k+1$. In principle, these values of $k_0$ and $\rho$ could be chosen before the induction process. Thus the sequence of functions described in (\ref{5.4}) has been constructed.

For any $x\in B(r_{k_0})$, there exists such a number $k$ that $x\in B(r_{kk_0})\setminus B(r_{(k+1)k_0})$, that is
$$
q^{-\lambda_{(k+1)k_0}}<|x|_K\le q^{-\lambda_{kk_0}}
$$
and
$$
|u(x)|\le \widehat{C}Mq^{-\rho \lambda_{kk_0}}.
$$
As we know, $\lambda_{kk_0}\ge \dfrac{\lambda_{(k+1)k_0}}{R_+^{k_0}}$, so that
$$
|u(x)|\le \widehat{C}Mq^{-\frac{\rho}{R_+^{k_0}} \lambda_{(k+1)k_0}}\le  \widehat{C}M|x|_K^{\frac{\rho}{R_+^{k_0}}},
$$
which means the H\"older continuity of the function $u$. $\qquad \blacksquare$

\bigskip
{\bf Examples.} 1) Let
$$
\Omega =\bigcup\limits_{k=1}^\infty S(r_k),\quad S(r_k)=\left\{ x\in K:\ |x|_K=q^{-\lambda_k}\right\}
$$
where $r_k=q^{-\lambda_k}$, and the sequence $\{ \lambda_k\}$ satisfies the inequalities (\ref{5.3}) with $R_->2$.

We have $0\in \partial \Omega$, $\Omega^c=\bigcap\limits_{k=1}^\infty S(r_k)^c=\{ x\in K:\ |x|_K\ne r_k \text{ for all $k=1,2,\ldots$} \}$,
$$
\left[ B(r_k)\setminus  B(r_{k+1})\right]\cap \Omega^c =\left\{ x\in K:\ r_{k+1}<|x|_K<r_k\right\},
$$
so that
$$
\operatorname{mes}\left[ B(r_k)\setminus  B(r_{k+1})\right]\cap \Omega^c=q^{-\lambda_k-1}-q^{-\lambda_{k+1}}=q^{-\lambda_k}\left( q^{-1}-q^{-(\lambda_{k+1}-\lambda_k)}\right).
$$

By our definition,
$$
\lambda_{k+1}-\lambda_k=\lambda_k\left( \frac{\lambda_{k+1}}{\lambda_k}-1\right)\ge (R_- -1)\lambda_k,
$$
and therefore
$$
q^{-1}-q^{-(\lambda_{k+1}-\lambda_k)}\ge q^{-1}-q^{-(R_- -1)\lambda_k}\ge q^{-1}-q^{-R_-+1}=q^{-1}\left( 1-q^{-(R_--2)}\right) \overset{\text{def}}{=} \nu >0,
$$
that is the condition (\ref{5.2}) is satisfied, and the solution on $\Omega$ is H\"older continuous at the origin.

\medskip
2) An obvious example of irregular behavior at the boundary is the punctured disk $\Omega =O\setminus \{ 0\}$ and the fundamental solution $u(x)=\dfrac{1-q^{-\alpha}}{1-q^{-\alpha -1}}|x|_K^{\alpha -1}$ (see \cite{K2001}, Section 2.2).

In this case, $\Omega^c=\{ x\in K:\ x=0 \text{ or $|x|_K>1$}\}$. Here the condition (\ref{5.2}) is of course violated.

\bigskip
-{\bf Ethical Approval}: The paper complies with ethical standards.  There are no conflicts of interests.

-{\bf Consent to Participate}:  Not applicable (there is a single author)

-{\bf Consent to Publish}: Not applicable (there is a single author)

-{\bf Authors Contributions}: Not applicable (there is a single author)

-{\bf Funding}: This work was funded in part under the State Budget of Ukraine, Program 6541230 "Support to Priority Trends of Scientific Research", and also under a fellowship of the Universities for Ukraine (U4U) program.

-{\bf Competing Interests}: There are no competing interests.

-{\bf Availability of data and materials}: Not applicable


\medskip
\begin{thebibliography}{999}
\bibitem{AKS}
S. Albeverio, A. Yu. Khrennikov, and V. M. Shelkovich, {\it Theory of $p$-Adic Distributions}, Cambridge University Press, 2010.
\bibitem{BCW}
A. Bendikov, W. Cygan and W. Woess, Oscillating heat kernels on ultrametric spaces, {\it J. Spectral Theory} {\bf 9} (2019), 195--226.
\bibitem{BGHH}
A. Bendikov, A. Grigor'yan, E. Hu and J. Hu, Heat kernels and nonlocal Dirichlet forms on ultrametric spaces, {\it  Ann. Sc. Norm. Super. Pisa} {\bf 22} (2021), 399–-461.
\bibitem{BGPW}
A. D. Bendikov, A. A. Grigor'yan, Ch. Pittet and W. Woess, Isotropic Markov semigroups on ultra-metric spaces, {\it Russian Math. Surveys}, {\bf 69} (2014), 589--680.
\bibitem{Bo}
N. Bourbaki, {\it General Topology: Chapters 1–4}, Springer, Berlin, 1989.
\bibitem{Br}
H. Brezis, {\it Functional Analysis, Sobolev Spaces and Partial Differential Equations}, Springer, New York, 2011.
\bibitem{BV}
C. Bucur and E. Valdinoci, {\it Nonlocal Diffusions and Applications}, Springer, Switzerland, 2016.
\bibitem{GT}
D. Gilbarg and N. S. Trudinger, {\it Elliptic Partial Differential Equations of Second Order}, Springer, Berlin, 1983.
\bibitem{GK1}
P. G\'orka and T. Kostrzewa, Sobolev spaces on metrizable groups, {\it Ann. Acad. Sci. Fenn. Math.} {\bf 40} (2015), 837--849.
\bibitem{GK2}
P. G\'orka and T. Kostrzewa, A second look of Sobolev spaces on metrizable groups, {\it Ann. Acad. Sci. Fenn. Math.} {\bf 45} (2020), 95--120.
\bibitem{GKR}
P. G\'orka, T. Kostrzewa and E. G. Reyes, Sobolev spaces on locally compact Abelian groups: compact embeddings and local spaces, {\it J. Function Spaces}, 2014, Article 404738, 6 p.
\bibitem{H}
S. Haran, Analytic potential theory over the p-adics, {\it Ann. Inst. Fourier} {\bf 43} (1993), 905--944.
\bibitem{Ka}
H. Kaneko, Besov space and trace theorem on a local field and its application, {\it Math. Nachr.} {\bf 285} (2012), 981--996.
\bibitem{KKZ}
A. Yu. Khrennikov, S. V. Kozyrev and W. A. Z\'u\~niga-Galindo, {\it Ultrametric Pseudodifferential Equations and Applications}, Cambridge University Press, 2018.
\bibitem{Kim}
Y.-C. Kim, A simple proof of the p-adic version of the Sobolev embedding theorem, {\it Commun. Korean Math. Soc.} {\bf 25} (2010), No. 1, 27--36.
\bibitem{K2001}
A. N. Kochubei, {\it Pseudo-Differential Equations and Stochastics
over Non-Archimedean Fields}, Marcel Dekker, New York, 2001.
\bibitem{K2014}
A. N. Kochubei, Radial solutions of non-Archimedean pseudodifferential equations, {\it Pacif. J. Math.} {\bf 269} (2014), 355--369.
\bibitem{K2018}
A. N. Kochubei, Linear and nonlinear heat equations on a p-adic ball, {\it Ukrainian Math. J.} {\bf 70} (2018), 217--231.
\bibitem{K2020}
A. N.  Kochubei, Non-Archimedean radial calculus: Volterra operator and Laplace transform, {\it Integr. Equat. Oper. Theory} {\bf 92: 44} (2020), 17 pp.
\bibitem{K2021}
A. N.  Kochubei, $L^p$ properties of non-Archimedean fractional differential operators, {\it J. Pseudodiff. Operators Appl.} {\bf 12}, No. 4 (2021), Article 56, 14 p.
\bibitem{LLF}
M. L. Lapidus, H\`ung L{\~u}\'{} and M. van Frankenhuijsen, Minkovski dimension and explicit tube formulas for p-adic fractal strings, {\it Fractal and Fractional}, 2018, 2, 26.
\bibitem{LZLH}
Y. Lian, K. Zhang, D. Li and G. Hong, Boundary H\"older regularity for elliptic equations, {\it J. Math. Pures Appl.} {\bf 143} (2020), 311--333.
\bibitem{LL}
E. Lindgren and P. Lindqvist, Fractional eigenvalues, {\it Calc. Var. Partial Diff. Equations}, {\bf 49} (2014), 795--826.
\bibitem{Ma}
V. G. Maz'ya, {\it Boundary Behavior of Solutions to Elliptic Equations in General Domains}, EMS, Z\"urich, 2018.
\bibitem{Po}
A. C. Ponce, {\it Elliptic PDEs, Measures and Capacities}, EMS, Z\"urich, 2015.
\bibitem{RO}
X. Ros-Oton, Nonlocal elliptic equations in bounded domains: a survey, {\it Publ. Mat.} {\bf 60} (2016), 3--26.
\bibitem{Sch}
W. H. Schikhof, {\it Ultrametric Calculus}, Cambridge University Press, Cambridge, 1984.
\bibitem{Se}
J.-P. Serre, {\it Lie Algebras and Lie Groups}, Benjamin, New York, 1965.
\bibitem{Se1}
J.-P. Serre, {\it Local Fields}, Springer, New York, 1979.
\bibitem{T}
M. H. Taibleson, {\it Fourier Analysis on Local Fields}, Princeton University Press, 1975.
\bibitem{Ta}
M. H. Taibleson, The existence of natural field structures for finite dimensional vector spaces over local fields, {\it Pacif. J. Math.}, {\bf 63} (1976), 545--551.
\bibitem{TJZ}
K. Teng, H. Jia and H. Zheng, Existence and multiplicity results for fractional differential inclusions with Dirichlet boundary conditions, {\it Appl. Math. Comput.} {\bf 220} (2013), 792--801.
\bibitem{TZ}
A. Torresblanca-Badillo and W. A. Z\`u\~{n}iga-Galindo, Non-Archimedean pseudodifferential operators and Feller semigroups, {\it p-Adic numbers, Ultrametric Anal. Appl.} {\bf 10}, No. 1 (2018). 57--73.
\bibitem{VVZ}
V. S. Vladimirov, I. V. Volovich and E. I. Zelenov, {\it $p$-Adic Analysis and
Mathematical Physics}, World Scientific, Singapore, 1994.
\bibitem{We}
A. Weil, {\it Basic Number Theory}, Springer, Berlin, 1967.
\bibitem{Z}
W. A. Z\`u\~{n}iga-Galindo, {\it  Pseudodifferential Equations over Non-Archimedean Spaces}, Lect.
Notes Math. Vol. 2174 (2016), XVI+175 p.
\bibitem{Z2022}
W. A. Z\`u\~{n}iga-Galindo, Ultrametric diffusion, rugged energy landscapes and transition networks. {\it Physica A} {\bf 597} (2022), Article 127221, 19 p.
\end{thebibliography}
\end{document}